\documentclass{amsart}
\usepackage{amsmath,amssymb,amscd,amsfonts}

\newtheorem{theorem}{Theorem}
\newcommand{\bt}{\begin{theorem}}
\newcommand{\et}{\end{theorem}}
\newtheorem{lemma}{Lemma}
\newcommand{\bl}{\begin{lemma}}
\newcommand{\el}{\end{lemma}}
\newtheorem{corollary}{Corollary}
\newcommand{\bc}{\begin{corollary}}
\newcommand{\ec}{\end{corollary}}
\newcommand{\beq}{\begin{equation}}
\newcommand{\eeq}{\end{equation}}
\newcommand{\benum}{\begin{enumerate}}
\newcommand{\eenum}{\end{enumerate}}
\newcommand{\N}{\ensuremath{ \mathbf N }}
\newcommand{\Z}{\ensuremath{\mathbf Z}}

\newcommand{\mcf}{\ensuremath{ \mathcal F}}

\DeclareMathOperator{\qand}{\quad\text{and}\quad}

\title{Sumsets contained in sets  of upper Banach density 1}

\author{Melvyn B. Nathanson}
\address{Lehman College (CUNY),Bronx, New York 10468}
\email{melvyn.nathanson@lehman.cuny.edu}

 \thanks{Supported in part by a grant from the PSC-CUNY Research Award Program.}

\date{\today}

\begin{document}

\begin{abstract}
Every set $A$ of positive integers with upper Banach density 1  
contains an infinite sequence of pairwise disjoint subsets 
$(B_i)_{i=1}^{\infty}$ such that $B_i$ has upper Banach density 1 for all $i \in \N$ 
and
$\sum_{i\in I} B_i  \subseteq A$ 
for every nonempty finite set $I$ of positive integers. 
\end{abstract}

\maketitle

\section{Upper Banach density}
Let  $\N$, $\N_0$, and \Z\ denote, respectively,  the sets of positive integers, nonnegative integers, and integers. 
Let $|S|$ denote the cardinality of the set $S$.
We define the \emph{interval of integers} 
\[
[x,y] = \{n \in \N : x \leq n \leq y\}.
\]

Let $A$ be a set of positive integers. 
Let $n \in \N$.  For all $u \in \N_0$, we have 
\[
|A \cap [u,u+n-1] | \in [0,n] 
\]
and so
\[
f_A(n) =  \max_{u \in \N_0}  \left| A \cap [u,u+n-1] \right|
\]
exists.  The \emph{upper Banach density} of $A$ is
\[
\delta(A) = \limsup_{n \rightarrow \infty} \frac{f_A(n)}{n}.
\]
Let  $n_1,n_2 \in \N$.  
There exists $u^*_1 \in \N_0$ such that, with $u^*_2 = u^*_1 + n_1$,  
\begin{align*}
f_A(n_1+n_2) & =  \left| A \cap [u^*_1, u^*_1+ n_1 + n_2 -1] \right| \\
& =  \left| A \cap [u^*_1, u^*_1+ n_1 -1] \right| +  \left| A \cap [u^*_1+n_1, u^*_1+ n_1 + n_2 -1] \right| \\
& =  \left| A \cap [u^*_1, u^*_1+ n_1 -1] \right| +  \left| A \cap [u^*_2, u^*_2 + n_2 -1] \right| \\
& \leq f_A(n_1)  + f_A(n_2). 
\end{align*}
It is well known, and proved in the Appendix, 
that this inequality implies that 
\[
\delta(A) = \lim_{n \rightarrow \infty} \frac{f_A(n)}{n} = \inf_{n \in \N} \frac{f_A(n)}{n}.
\]

\section{An Erd\H os sumset conjecture}
About 40 years ago, Erd\H os conjectured that if $A$ is a set of positive integers of positive 
upper Banach density, then there exist infinite sets $B$ and $C$ 
of  positive integers such that $B+C \subseteq A$.  
This conjecture has not yet been verified or disproved.  

The \emph{translation}\index{translation} of the set $X$ by $t$ is the set 
\[
X+t = \{x + t : x\in X\}.
\]
Let $B$ and $C$ be sets of integers.  
For every integer $t$, if $B' = B+t$ and $C' = C - t$,   
then 
\[
B' + C' = (B+t) + (C-t) = B+C.
\]
In particular, if $C$ is bounded below and $t = \min(C)$, 
then $0 = \min(C')$ and $B' \subseteq B' + C'$.  
It follows that if $B$ and $C$ are infinite sets such that $B+C\subseteq A$, 
then, by translation, there exist  infinite sets $B'$ and $C'$ such that 
$B' \subseteq A$ and $B' +C' \subseteq A$.  

However, a set $A$ with positive upper Banach density does not necessarily 
contain infinite subsets $B$ and $C$ with $B+C \subseteq A$.  
For example, let $A$ be any set of odd numbers.  
For all sets $B$ and $C$  of odd numbers, the sumset  $B+C$ is a set of even numbers, 
and so $A \cap (B+C) = \emptyset$.  Of course, in this example we have 
$B+C \subseteq A+1$.  

In this note we prove that if $A$ is a set of positive integers 
with upper Banach density $\delta(A) = 1$, then 
for every $h \geq 2$  there exist pairwise disjoint subsets $B_1,\ldots, B_h$ of $A$ 
such that  $\delta(B_i) = 1$ for all $i=1,\ldots, h$ and 
\[
B_1 + \cdots + B_h \subseteq A.
\]
Indeed, Theorem~\ref{sumsets:theorem:I-sequence} states an even stronger result.

There are sets $A$ of upper Banach density 1 for which no infinite subset $B$ of $A$ satisfies 
$2B \subseteq A + t$ for any integer $t$.  A simple example is 
\[
A = \bigcup_{i=1}^{\infty} \left[ 4^{i}, 4^{i}  + i - 1\right].
\]
The set $A$ is the union of the infinite sequence of pairwise disjoint intervals 
\[
A_i = \left[ 4^{i}, 4^{i} + i - 1\right].
\]
Let $t \in \N_0$.  There exists $i_0(t)$ such that $4^i-i > t$ for all $i \geq i_0(t)$.  
If $b_i \in A_i$ for some $i \geq i_0(t)$, then 
\[
4^i +i +  t < 2\cdot 4^{i} \leq 2b_i <  2\cdot 4^{i} + 2i  < 4^{i+1} - 2t \leq 4^{i+1} - t
\]
and so $2b_i \notin 2A\pm t$.  
If $B$ is an infinite subset of $A$, then for infinitely many $i$ there exist integers 
$b_i \in B \cap A_i$, and so $2B \not\subseteq A+t$ for all $t \in \Z$.

There are very few results about the Erd\H os conjecture.  
In 1980, Nathanson~\cite{nath80d} proved that if $\delta(A) > 0$, then for every $n$ 
there is a finite set $C$ with $|C| = n$ and a subset $B$ of $A$ with 
$\delta(B) > 0$ such that  $B+C\subseteq A$.
In 2015, Di Nasso, Goldbring, Jin, Leth, Lupini, and Mahlburg~\cite{dina-gold15} 
used nonstandard analysis to prove that 
the Erd\H os conjecture is true for sets $A$ with upper Banach density $\delta(A) > 1/2$.  
They also proved that if $\delta(A) > 0$, then there exist infinite sets $B$ and $C$ 
and an integer $t$ such that 
\[
B+C \subseteq A \cup (A+t).
\]
It would be of interest to have purely combinatorial proofs of the 
results of Di Nasso, \emph{et al.}

For related work, see Di Nasso~\cite{dina14a,dina14b}, 
Gromov~\cite{grom15}, Hegyv\' ari~\cite{hegy99,hegy08},  
Hindman~\cite{hind79a}, and Jin~\cite{jin04}.

\section{Results}
The following result is well known.  

\bl        \label{sumset:lemma:Banach}
A  set of positive integers has upper Banach density 1 
if and only if, for every $d$,  it contains 
infinitely many pairwise disjoint intervals of $d$ consecutive integers.   
\el

\begin{proof}
Let $A$ be a set of positive integers.  If, for every positive integer $d$, 
the set $A$ contains an interval of $d$ consecutive integers, then
\[
\max_{u \in \N_0} \left( \frac{ |A \cap [u,u+d-1] |}{d} \right) = 1
\]
and so
\[
\delta(A) = \lim_{d \rightarrow \infty} \max_{u \in \N_0} \left( \frac{ |A \cap [u,u+d-1] |}{d} \right) = 1.
\]

Suppose that, for some integer $d \geq 2$, the set $A$ contains no interval 
of $d$ consecutive integers.   
For every $u \in \N_0$, we consider the interval $I_{u,n} = [u, u+n-1]$. 
By the division algorithm, there are integers $q$ and $r$ with $0 \leq r < d$ 
such that 
\[
|I_{u,n}| = n = qd + r
\]
and  
\[
q = \frac{n-r}{d} > \frac{n}{d} - 1.
\]
For $j=1,\ldots, q$,  the intervals of integers  
\[
I_{u,n}^{(j)} = [u+(j-1)d, u + jd - 1 ]
\]
and 
\[
I_{u,n}^{(q+1)} = [u+qd, u+n-1]
\]
are pairwise disjoint subsets of $I_{u,n}$ such that 
\[
I_{u,n}  = \bigcup_{j=1}^{q+1}I_{u,n}^{(j)}. 
\] 
We have 
\[
A \cap I_{u,n}  = \bigcup_{j=1}^{q+1} ( A \cap I_{u,n}^{(j)} )
\]
If  $A$ contains no interval of $d$ consecutive integers, 
then, for all $j \in [1,q]$, at least one element of the interval $I_{u,n}^{(j)}$ is not 
an element of $A$, and so 
\[
| A \cap I_{u,n}^{(j)}|  \leq |I_{u,n}^{(j)}| - 1.  
\]
It follows that  
\begin{align*}
| A \cap  I_{u,n}|  & = \sum_{j=1}^{q+1} \left| A \cap  I_{u,n}^{(j)} \right| 
 \leq  \sum_{j=1}^q \left(  \left| I_{u,n}^{(j)}\right| - 1 \right)  + \left| I_{u,n}^{(q+1)}\right| \\
&  =  \sum_{j=1}^{q+1}  \left| I_{u,n}^{(j)}\right|  - q =  |I_{u,n}| - q  = n - q \\ 
& <   n -  \frac{n}{d} +1 = \left(1 - \frac{1}{d}\right)n + 1. 
\end{align*}
Dividing by $n = |I_{u,n}|$, we obtain 
\[
\max_{u \in \N_0} \frac{ | A \cap I_{u,n}|}{n} \leq 1 - \frac{1}{d} + \frac{1}{n}.
\]
and so 
\[
\delta(A) = \lim_{n \rightarrow \infty} \max_{u \in \N_0} \frac{ | A \cap  I_{u,n} |}{n}  
\leq 1 - \frac{1}{d} < 1
\]
which is absurd.  Therefore, $A$ contains an interval 
of $d$ consecutive integers for every $d \in \N$.  

To prove that $A$ contains infinitely many intervals of size $d$, it suffices to prove 
that if $[u,u+d-1]\subseteq A$, then $[v,v+d-1] \subseteq A$ for some $v \geq u+d$.
Let $d' = u+2d$.  There exists $u' \in \N$ such that 
\[
[u',u'+d'-1] = [u', u' + u + 2d-1] \subseteq A.
\]
Choosing $v = u'+u+d$, we have $v \geq u+d$ and 
\[
[v,v+d-1] \subseteq [u', u' + u + 2d-1] \subseteq A.
\] 
This completes the proof.  
\end{proof}

Let $\mcf(S)$ denote the set of all finite subsets of the set $S$, 
and let $\mcf^*(S)$ denote the set of all nonempty finite subsets of $S$.
We have the fundamental binomial  identity 
\beq              \label{sumsets:FiniteSubsets}
\mcf^*([1,n+1] ) = \mcf^*([1,n])  \cup \left\{ \{n+1\}\cup J:J\in \mcf(  [1,n]  )  \right\}.
\eeq

\bt                                \label{sumsets:theorem:b-sequence}
Let $A$ be a set of positive integers that  has upper Banach density 1.   
For every sequence $(\ell_j)_{j=1}^{\infty}$ of positive integers, 
there exists a sequence $(b_j)_{j=1}^{\infty}$ of positive integers 
such that 
\[
b_{j+1} \geq  b_{j} + \ell_j 
\]
for all $j \in \N$, and 
\[
 \sum_{j \in J} [b_{j}, b_{j} + \ell_j -1] \subseteq A
\]
for all $J  \in \mcf^*(\N)$.
\et

\begin{proof}
We shall construct the sequence $(b_j)_{j=1}^{\infty}$ by induction.  
For $n=1$, choose $b_1 \in A$ such that $[b_1,b_1+\ell_1 -1] \subseteq A$. 

Suppose that  $(b_j)_{j=1}^n$ is a finite sequence of positive integers 
such that $b_{j+1} \geq b_{j} + \ell_j$ for $j \in [1,n-1]$ and 
\beq              \label{sumsets:FiniteSubsets-2}
\sum_{j \in J}  [b_{j}, b_{j} + \ell_j -1]  \subseteq A
\eeq       
for all $J  \in \mcf^*([1,n])$.
By Lemma~\ref{sumset:lemma:Banach}, 
there exists $b_{n+1} \in A$ such that 
\[
b_{n+1} \geq  b_n + \ell_n 
\]
and
\[
\left[  b_{n+1}, \sum_{j=1}^{n+1} (b_j + \ell_j ) -1\right] \subseteq A. 
\]
It follows that 
\[
 \left[  b_{n+1}, b_{n+1} + \ell_{n+1}-1\right] \subseteq A. 
\]

Let $J  \in \mcf([1,n])$.  If 
\begin{align*}
a  & \in   \sum_{j \in \{n+1\} \cup J} [b_{j}, b_{j} + \ell_j -1] \\
& =  \left[  b_{n+1}, b_{n+1} + \ell_{n+1}-1\right]  +  \sum_{j \in J} [b_{j}, b_{j} + \ell_j -1]
\end{align*}
then 
\begin{align*}
b_{n+1}  \leq a & \leq \left(  b_{n+1} + \ell_{n+1}-1 \right) 
+  \sum_{j \in J} \left( b_{j} + \ell_j -1 \right)  \\
& \leq \sum_{j=1}^{n+1} (b_j + \ell_j ) -1
\end{align*}
and so $a \in A$ and 
\beq              \label{sumsets:FiniteSubsets-3}
 \sum_{j \in \{n+1\} \cup J} [b_{j}, b_{j} + \ell_j -1] 
 \subseteq \left[  b_{n+1}, \sum_{j=1}^{n+1} (b_j + \ell_j ) -1\right] \subseteq A. 
\eeq
Relations~\eqref{sumsets:FiniteSubsets},~\eqref{sumsets:FiniteSubsets-2}, 
and~\eqref{sumsets:FiniteSubsets-3} imply that 
\[
\sum_{j \in J}  [b_{j}, b_{j} + \ell_j -1]  \subseteq A
\]      
for all $J  \in \mcf^*([1,n+1])$.
This completes the induction.  
\end{proof}

\bt                                               \label{sumsets:theorem:I-sequence}
Every set $A$ of positive integers that  has upper Banach density 1  
contains an infinite sequence of pairwise disjoint subsets 
$(B_i)_{i=1}^{\infty}$ such that $B_i$ has upper Banach density 1 for all $i \in \N$ 
and
\[
\sum_{i\in I} B_i  \subseteq A
\] 
for all $I \in \mcf^*(\N)$.  
\et

\begin{proof}
Let $(\ell_j)_{j=1}^{\infty}$ be a sequence of positive integers such that 
$\lim_{j\rightarrow \infty} \ell_j = \infty$, and let $(b_j)_{j=1}^{\infty}$ 
be a sequence of positive integers that satisfies 
Theorem~\ref{sumsets:theorem:b-sequence}.  
(For simplicity, we can let $\ell_j = j$ for all $j$.)
Let $(X_i)_{i =1}^{\infty}$ be a sequence of  infinite sets of positive integers 
that are pairwise disjoint.  
For $i \in \N$, let 
\[
B_i = \bigcup_{j\in X_i} [b_{j}, b_{j} + \ell_j -1].  
\]
The set $B_i$ contains  intervals of $\ell_j$ consecutive integers for infinitely 
many $\ell_j$, and so $B_i$ has upper Banach density 1.  

Let $I \in \mcf^*(\N)$.  If 
\[
a \in \sum_{i\in I} B_i  \subseteq A
\] 
then for each $i \in I$ there exists $a_i \in B_i$ such that $a = \sum_{i\in I} a_i$.  
If $a_i \in B_i$, then there exists $j_i \in X_i$ such that 
\[
x_i \in \left[ b_{j_i}, b_{j_i} + \ell_{j_i}  - 1 \right].
\]
We have $J = \{ j_i: i \in I\} \in \mcf^*(\N)$ and 
\[
a \in \sum_{j_i \in J} \left[ b_{j_i}, b_{j_i} + \ell_{j_i}  - 1 \right] \subseteq A.
\]
This completes the proof.  
\end{proof}

\bt                                               \label{sumsets:theorem:I-sequence-AP}
Let $A$ be a set of integers that contains arbitrarily long finite arithmetic progressions 
with bounded differences.  
There exist positive integers $m$ and $r$, and an infinite sequence 
of pairwise disjoint sets 
$(B_i)_{i=1}^{\infty}$ such that $B_i$ has upper Banach density 1 for all $i \in \N$ 
and
\[
m\ast \sum_{i\in I}  B_i +r \subseteq A
\] 
for all $I \in \mcf^*(\N)$.  
\et

\begin{proof}
If the differences in the infinite set of finite arithmetic progressions 
contained in $A$ are bounded by $m_0$, then 
there exists a difference $m \leq m_0$ that occurs infinitely often.  
It follows that there are 
arbitrarily long finite arithmetic progressions with  difference $m$.
Because there are only finitely many congruence classes modulo $m$, 
there exists a congruence class $r \pmod{m}$ such that $A$ contains 
arbitrarily long sequences of consecutive integers in the the congruence class 
$r\pmod{m}$.  Thus, there exists an infinite set $A'$ such that 
\[
m\ast A' + r \subseteq A
\]
and $A'$ contains arbitrarily long sequences of consecutive integers.  
Equivalently, $A'$ has Banach density 1.
By Theorem~\ref{sumsets:theorem:I-sequence}, the sequence $A'$ 
contains an infinite sequence of pairwise disjoint subsets 
$(B_i)_{i=1}^{\infty}$ such that $B_i$ has upper Banach density 1 for all $i \in \N$ 
and
\[
\sum_{i\in I} B_i  \subseteq A'
\] 
for all $I \in \mcf^*(\N)$.  
It follows that 
\[
m\ast \sum_{i\in I} B_i +r \subseteq m\ast A' + r \subseteq A
\] 
for all $I \in \mcf^*(\N)$.  
This completes the proof.  
\end{proof}

\appendix
\section{Subadditivity and limits} \label{sumset:appendix}

A real-valued arithmetic function $f$ is \emph{subadditive}\index{subadditive} if 
\beq   \label{sumset:subadditive}
f(n_1+n_2) \leq f(n_1)  + f(n_2) 
\eeq
for all $n_1,n_2 \in \N$.

The following result is sometimes called \emph{Fekete's lemma}.

\bl         \label{sumset:lemma:subadditiveLimit}
If $f$ is a subadditive arithmetic function, 
then $\lim_{n\rightarrow \infty} f(n)/n$ 
exists, and 
\[
\lim_{n\rightarrow \infty} \frac{f(n)}{n} = \inf_{n\in \N} \frac{f(n)}{n}.
\]
\el

\begin{proof}
It follows by induction from inequality~\eqref{sumset:subadditive} 
that 
\[
f(n_1+\cdots + n_q) \leq f(n_1) + \cdots + f(n_q)
\]
for all $n_1,\ldots, n_q \in \N$.  Let $f(0) = 0$.  Fix a positive integer $d$.  
For all $q,r \in \N_0$,  we have
\[
f(qd+r) \leq qf(d)+ f(r).
\]
By the division algorithm, every  nonnegative 
integer $n$ can be represented uniquely in the form $n = qd+r$, 
where $q \in \N_0$ and $r \in [0,d-1]$.  Therefore,
\[
\frac{f(n)}{n} = \frac{f(qd+r)}{n} \leq  \frac{qf(d)}{qd} + \frac{f(r)}{n} 
=  \frac{f(d)}{d} + \frac{f(r)}{n}. 
\]
Because the set $\{f(r):r\in [0,d-1]\}$ is bounded, it follows that  
\[
\limsup_{n\rightarrow \infty} \frac{f(n)}{n}  \leq  
\limsup_{n\rightarrow \infty} \left(  \frac{f(d)}{d} + \frac{f(r)}{n} \right) 
= \frac{f(d)}{d} 
\]
for all $d \in \N$, and so 
\[
\limsup_{n\rightarrow \infty} \frac{f(n)}{n}  \leq   \inf_{d\in \N} \frac{f(d)}{d} 
\leq \liminf_{d \rightarrow \infty} \frac{f(d)}{d} =  \liminf_{n\rightarrow \infty} \frac{f(n)}{n}.  
\]
This completes the proof.  
\end{proof}

\def\cprime{$'$} \def\cprime{$'$} \def\cprime{$'$}
\providecommand{\bysame}{\leavevmode\hbox to3em{\hrulefill}\thinspace}
\providecommand{\MR}{\relax\ifhmode\unskip\space\fi MR }
\providecommand{\MRhref}[2]{%
  \href{http://www.ams.org/mathscinet-getitem?mr=#1}{#2}
}
\providecommand{\href}[2]{#2}

\end{document}